\documentclass[preprint,12pt]{elsarticle}
\usepackage{amsfonts, amsmath,latexsym,graphics,amsthm,amssymb} 
\newcounter{mascotsection}
\newcounter{mascotsubsection}[mascotsection]

\newcommand{\mascotsection}[1]{%
	\stepcounter{mascotsection}
	\section{\themascotsection . #1}
}

\usepackage{epsfig}

\newtheorem{theorem}{Theorem}
\newtheorem{dhef}{Definition}

\def\erre{{\rm I}\!{\rm R}}

\journal{Proceedings of MASCOT2015 Workshop, Rome \ }
\begin{document}

\begin{frontmatter}
\title{A numerical study of a two-layer model \\ for the growth of granular matter in a silo}
\author[rm]{S. Finzi Vita}
\address[rm]{Dipartimento di Matematica, Sapienza Universit\`a di Roma\\ P.le Aldo Moro, 5 - 00185 Roma, Italy\\ Email address: finzi@mat.uniroma1.it}

\begin{abstract}
\noindent The problem of filling a silo of given bounded cross-section with granular matter can be described by the two-layer model of Hadeler and Kuttler \cite{HK99}. In this paper we discuss how similarity quasi-static solutions for this  model can be  numerically characterized by the direct finite element solution of a semidefinite elliptic Neumann problem. We also discuss a finite difference scheme for the dynamical model through which we can show that the growing profiles of the heaps in the silo evolve in finite time towards such similarity solutions.
\end{abstract}
\end{frontmatter}
%
\medskip
\mascotsection{Introduction}
In recent times many models have been proposed to study the dynamics of granular materials like sand in different practical situations (see e.g. \cite{AT} for an overview). In many applications one has to store this kind of materials in view of their later use, so that filling and emptying a container are crucial processes. Granular materials adapt their shape to the container (like a fluid does), but in general the free surface of a heap strongly depends on its formation process, for example on intensity and dislocation of the source.  Moreover, the pressure on the bottom of the structure does not grow linearly with the height of the pile, since part of it is released against the walls through arcs of grains, a fact that can even produce silos explosion and collapse.

 In this paper we deal with the simple problem of pouring at low intensity granular matter into a silo of given cross-section $\Omega\subset \erre^2$: it is known from the experiments (see e.g. \cite{GH99}) that if the source is independent of time,  the free surface of the growing heap evolves towards a well-defined profile which then retains its shape while growing with a constant velocity. 
 
The two-layer model of Hadeler and Kuttler \cite{HK99}, which basically describes the formation of sandpiles over an open bounded plane table, is a system of two partial differential equations for a standing layer $u$ and a small rolling layer $v$ of grains running down the slope. It can be adapted  in a natural way to the case of the silo problem (see again \cite{HK99} and more specifically \cite{HK99b}) by adding a suitable boundary condition on the silo walls. If $f$ denotes the vertical source of material and $T$ the final time, then the model takes the form
\begin{equation}\label{hk}
 \left\{\begin{array}{ll}
         \displaystyle v_t = \beta\nabla \cdot (v\: \nabla u) - \gamma( \alpha-\vert \nabla u\vert)\: v + f &\hbox{ in } \Omega\times (0,T) \\  \\
\displaystyle u_t=  \gamma( \alpha-\vert \nabla u\vert)\: v &\hbox{ in } \Omega\times (0,T) \\ \\
u(x,0)=u_0(x)\ ,\quad  v(x,0)=0 &\hbox{ in } \Omega\\ \\
\displaystyle\frac{\partial u}{\partial n}=0  &\hbox{ on } \partial\Omega\times (0,T). 
        \end{array}
\right.
\end{equation}
The nonlinear term which appears in both the equations with opposite sign expresses the exchange term between the two layers, $\alpha$ being the maximal ({\it critical}) slope that the material can support without flowing down, $\beta$ and $\gamma$ respectively the mobility and the collision rate parameters. The boundary condition comes from the total mass conservation law (when $f=0$), which suggests
\begin{equation}\label{bc}
v\frac{\partial u}{\partial n} = 0 \hbox{ on } \partial\Omega \ ;
\end{equation}
but the first equation in (\ref{hk}) is an advection equation for $v$ in the direction of $(-\nabla u)$, so  that  boundary conditions  cannot be imposed for the outgoing direction of its flow and condition (\ref{bc}) reduces to a pure homogeneous Neumann condition on $u$. Existence and uniqueness results for the solution of a system similar to (\ref{hk}) under general assumptions on the data has been recently discussed in \cite{DPJ}. 
If the source is constant in time the free profile is expected to evolve towards a {\it similarity solution}, according to the following definition.

\begin{dhef} We call a pair of functions $(U(x),V(x))$ a similarity solution of system (\ref{hk}) in $\Omega$ if there exists a positive constant $c$ such that the functions
\begin{equation}\label{sim} 
u(x,t) = U(x) + ct\ , \quad v(x,t)= V(x) 
\end{equation}
solve the system (disregarding the initial condition).
\end{dhef}

\medskip It can be considered as a sort of equilibrium for the model: the rolling layer is constant in time, while the free surface keeps growing by a rigid traslation of its shape at a constant rate $c$.

In the toy one-dimensional (1D) case for the cross-section, these quasi-stationary profiles $(U,V)$ can be  expressed by closed integral formulas (see \cite{HK99} and next section) in terms of the source and of the other problem parameters. In two dimensions (2D) it can be proved that similarity solutions exist, but their expressions are known only in special cases. That is why in Section 3 we discuss a finite element (FE) characterization of such solutions in the general case.

Finite difference (FD) numerical schemes for the model of Hadeler and Kuttler have been studied in \cite{FFV06} and \cite{FFV08} in the case of growing sandpiles on a bounded open table, and in \cite{CFV} on a table partially bounded by vertical walls. In Section~4 we will adapt such schemes to the present problem of silos in order to show through the experiments of Section 5 that the growing heaps generated by the evolving model perfectly match the similarity solutions.

\mascotsection{Characterization of similarity solutions}

We recall the basic theorem of existence for similarity solutions in the case of a constant in time source term.
\begin{theorem} (\cite{HK99}) Assume $f=f(x)$; then there exists a similarity solution $(U,V)$ for problem (\ref{hk}) (in the sense of Definition 1), with $U$ unique up to an additive constant and
\begin{equation}\label{vel}
c= \frac{1}{|\Omega|} \int_\Omega f\ dx .
\end{equation}
\end{theorem}
The main idea in the proof is to consider the basic properties of the flux function $w=v\nabla u$. If one looks for it in the form of a gradient ($w=\nabla \psi$), then its potential $\psi$ should solve the  semidefinite Neumann problem for the Laplacian
\begin{equation}\label{Neum}
       -\Delta \psi =  g \hbox{ in } \Omega, \qquad
      \frac{\partial \psi}{\partial n}  = 0 \hbox{ on } \partial\Omega,
\end{equation}
with $g=(f-c)/\beta$, which has a solution (unique up to an additive constant) due to the zero-mean property of  $g$. Then we can  derive $w$ from $\psi$, and also deduce:
\begin{equation}\label{vfromw} V = c +V|\nabla U| = \frac{1}{\gamma\alpha |\Omega|} \int_\Omega f\ dx + \frac 1 \alpha |w|,\qquad \nabla U= \frac{w}{V},\end{equation}
so that $U$ can be determined up to an additive constant.

Formula (\ref{vel}) says that the growth velocity of the similarity profile coincides with the average precipitation, that is with the mean value of the source intensity and is independent from the other parameters. Theorem 1 does not say anything about uniqueness: in principle other solution pairs $(U,V)$ could exist such that $V\nabla U$ is not a gradient. Anyway, numerical experiments of Section 5 show that the solutions given by the previous theorem are the only significant (physical) ones, since the evolving profiles tend asymptotically to them.\\
In 1D the previous result yields explicit expressions for similarity solutions. If for example $\Omega$ coincides with the interval $(0,L)$, one finds (see \cite{HK99} for details):
\begin{equation}\label{sim1}
 V(x)=\frac 1 {\gamma\alpha L}\int_0^L f(y)dy + \frac 1 {\alpha\beta}|G(x)|, \quad
 U_x(x)=\alpha \frac {G(x)} {\frac \beta {\gamma L}\int_0^L f(y)dy + |G(x)|},
\end{equation} 
where 
$$ G(x)=\frac x L \int_0^L f(y)dy - \int_0^x f(y)dy . $$
Such expressions give several informations about solutions:
\begin{itemize}
\item if $f(x)\equiv k\in\erre^+$ for any $x\in (0,L)$, then $G(x)\equiv 0$,  $V(x)=k/(\gamma\alpha)$ and $U_x \equiv 0$, that is the free surface grows remaining flat, as expected;
\item if the source is not identically zero, then $V(x)>0$ everywhere, even at the boundary, confirming what already stated about condition (\ref{bc});
\item $|U_x| \le \alpha$, that is the standing layer never exceeds the critical slope; 
\item the rolling layer thickness is directly proportional to the source intensity and 
inversely proportional to $\alpha$.
\end{itemize}
In higher dimensions explicit formulas for similarity solutions cannot be deduced in general, and we will see in the next section how to detect them numerically. Here we just report the two special cases of a central point source, for 1D and 2D radial cross-sections respectively, where similarity solutions can be explicitly computed.

\begin{figure}[h] 
\center{\fbox{\psfig{file=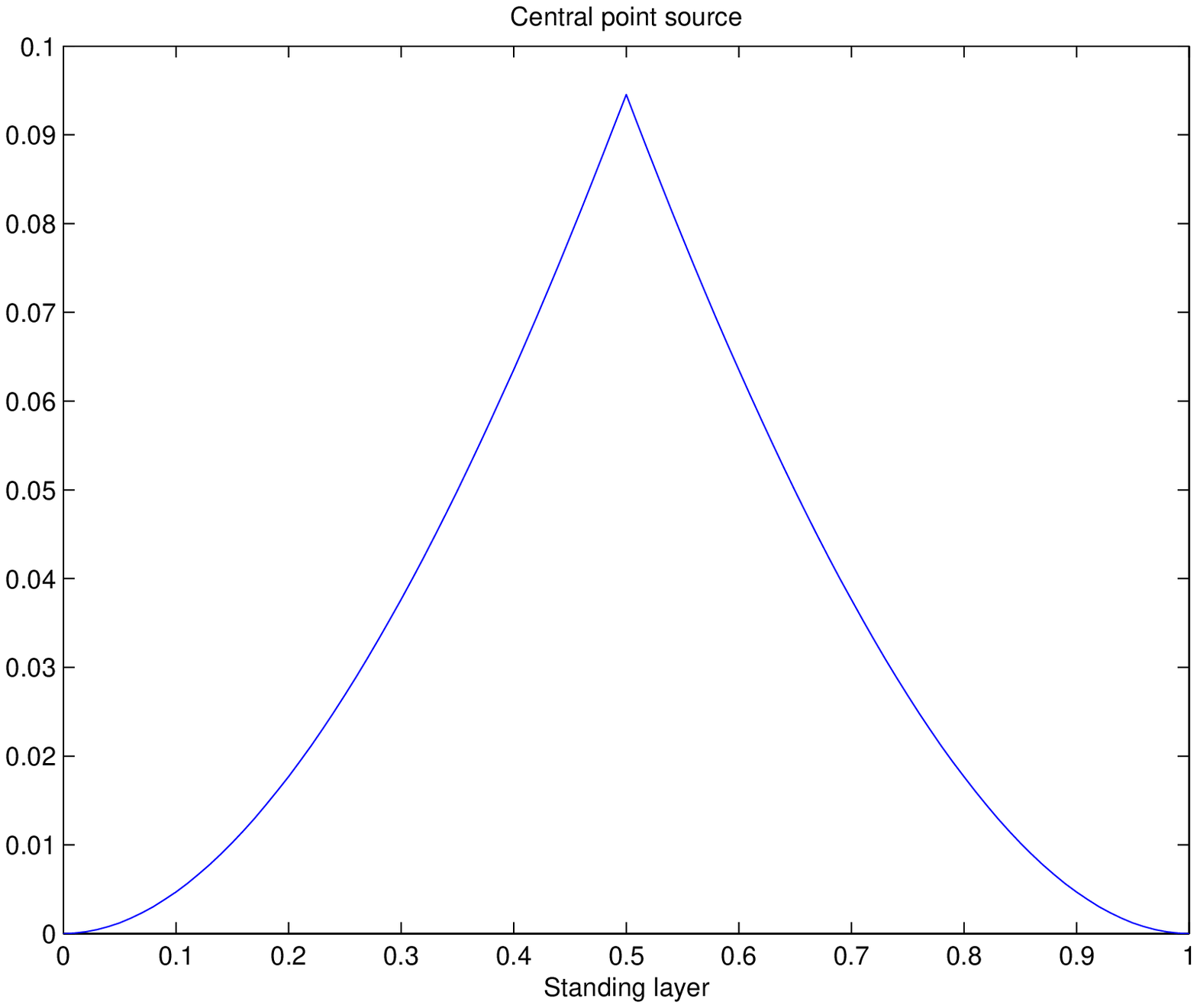,height=1.1in}}\hspace{.2in}
\fbox{\psfig{file=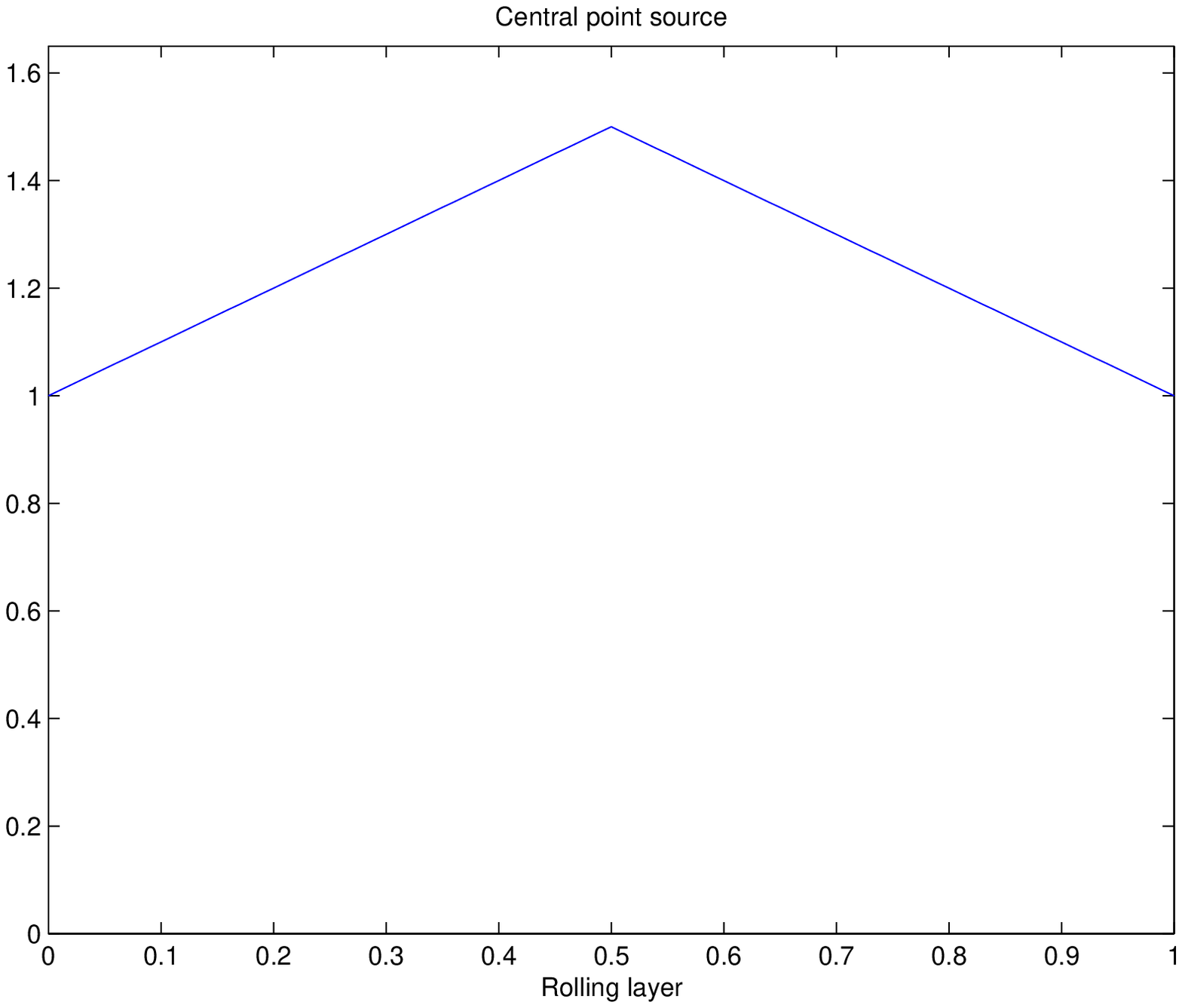,height=1.1in}}}
\caption{Similarity solutions $U$ and $V$ for 1D silo with central point  source.}
\label{sol1ddelta}
\end{figure}

\medskip\noindent{\bf Example 1.}(\cite{HK99}) Assume $\Omega=(0,L)$ and $f=\delta_{L/2}$ (where $\delta_z$ denotes the usual Dirac  function centered in $z$), that is there is a point source placed over the middle of the silo. Similarity solutions then take the form (see Figure \ref{sol1ddelta}):
\begin{equation}\label{dirac1_v}
V(x)=\frac 1 {\gamma\alpha L} + \frac 1 {\alpha\beta L}\min \{x,L-x\},
\end{equation} 
\begin{equation}\label{dirac1_u}
 U(x)=\left\{\begin{array}{ll}
     \displaystyle  \frac {\alpha\gamma} \beta (x-\log(1+x)) &\hbox{ if } x\le L/2\\  \\
    \displaystyle  \frac{\alpha\gamma} \beta (L-x-\log(1+L-x)) &\hbox{ if } x>L/2\ .
        \end{array}\right.
\end{equation}

\begin{figure}[t] 
\center{\fbox{\psfig{file=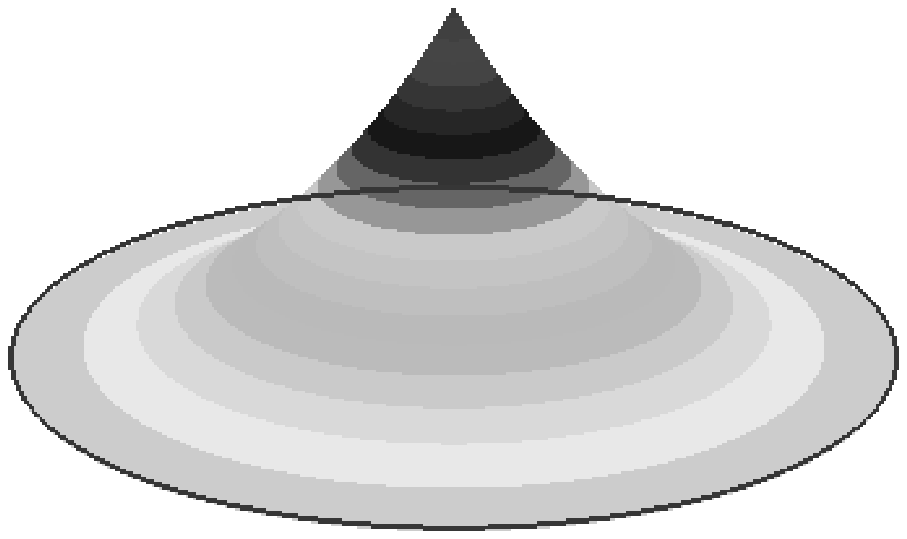,height=1.3in}}\hspace{.2in}
\fbox{\psfig{file=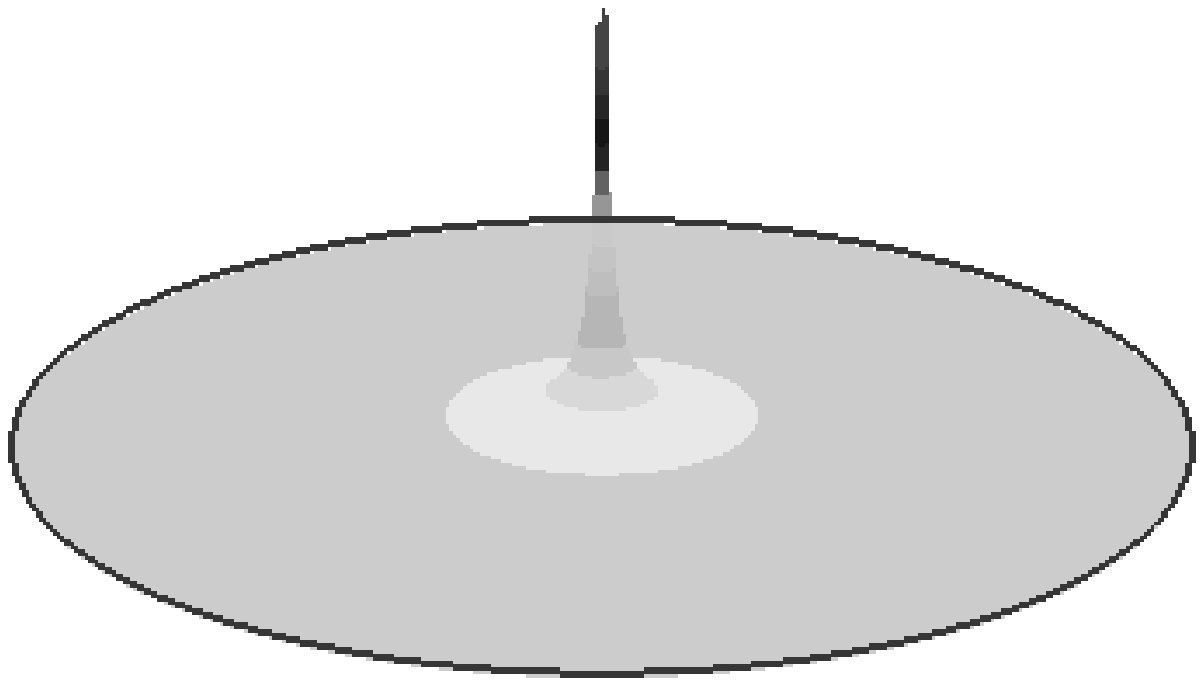,height=1.3in}}}
\caption{Similarity solutions $U$ and $V$ for a cylindrical silo with central point source.}
\label{sol3d}
\end{figure}

\noindent{\bf Example 2.}(\cite{HK99b}) Let $\Omega\subset\erre^2$ be the ball of radius $R$ centered at the origin, and $f=\delta_{(0,0)}$  a point source over its center;  then similarity solutions are radial functions, and radial symmetry arguments yields (see Figure \ref{sol3d}):
\begin{equation}\label{sim1R}
V(r)=\frac c {\gamma\alpha}(1+ \frac \gamma {2\beta r}(R^2-r^2)), \quad
U_r(r)=-\alpha \frac {R^2-r^2} {R^2-r^2+2\beta r/\gamma}.
\end{equation} 
Previous examples show that the typical profile of a growing symmetric heap of grains in a silo under the effect of a central source is different from the classical conical shape which would emerge on an open table without walls. In both dimensions the standing layer now assumes a strictly convex (logarithmic) profile. The maximal slope of the pile changes instead with the dimension: in 2D it is reached right in the center, and coincides with the critical slope $\alpha$, whereas in 1D it depends on the parameters $\alpha, \beta, \gamma$ and on the size $L$ of the container. In particular, when the ratio $\beta/\gamma <<1$ (grains roll slowly and are easily trapped) the slope of the pile always remains very close to the critical angle $\alpha$. On the contrary, when $\beta/\gamma$ is large the grains move very fast from the beginning, and larger variations of the slope can emerge.

For what concerns the rolling layer in 2D, Figure \ref{sol3d} shows the emergence of a singularity  in the center, in accordance with the fact that its expression comes from the solution of a potential problem with a central Dirac source.

\mascotsection{Approximation of similarity solutions}

For the sake of simplicity from now on we assume $\alpha=\beta=\gamma=1$.

The proof of Theorem 1 in the previous section shows that in order to  characterize the similarity profiles one needs to solve in $\Omega$ the elliptic Neumann problem (\ref{Neum}). From a numerical point of view, this can be done for example by using a finite element approach.  If  $T_h$ denotes a regular triangulation of $\Omega$ of size $h>0$, and  $V_h\subset H^1(\Omega)$ and $W_h\subset L^2(\Omega)$ are the finite element spaces of respectively  piecewise linear and piecewise constant functions on $T_h$, Galerkin method requires to solve the discrete variational problem\begin{equation}\label{discr_neum}
\psi_h\in V_h, \quad \int_{\Omega} \nabla \psi_h \cdot \nabla \phi \ dx = \int_{\Omega} g \phi \ dx\ , \quad \forall \phi\in V_h.
\end{equation}
However, if $\Omega$ is not a polygonal domain it has to be replaced in (\ref{discr_neum}) by a suitable set $\Omega_h$ defined as the union of the triangular elements of $T_h$; $\Omega_h$ will be closed to $\Omega$, but in general the right-hand side $g$  will not retain its zero-mean property on $\Omega_h$, and the discrete semidefinite problem will not be solvable at all. A way to overcome this difficulty is to replace $g$ in (\ref{discr_neum}) by the function $ g_h= f-c_h$, with $c_h=\frac 1 {|\Omega_h|}\int_{\Omega_h} f dx\ ;$ by construction, $ \lim_{h\to 0} g_h=g$, and 
$\int_{\Omega_h} g_h dx=0$ for any $h$,
so that (\ref{discr_neum}) becomes solvable (see \cite{CF88} for details). Then by definition
$w_h=\nabla \psi_h$ will be a piecewise constant vector on $T_h$, and its norm an element of the discrete space $W_h$. Hence, from (\ref{vfromw}),
\begin{equation}
v_h=c_h+|w_h|=c_h + |\nabla \psi_h|\ .
\end{equation}
Now, since $w_h= v_h \nabla u_h$, the gradient of $u_h$ on any triangle $\tau_k\in T_h$ is given by
\begin{equation}
z_k= \nabla u_h |_{\tau_k}=\frac {w_h}{v_h}|_{\tau_k}.
\end{equation}
It remains to compute $u_h\in V_h$. In 1D its value at any node $x_i\in (0,L)$ can be determined by direct  integration of the piecewise constant function $(u_h)_x$ from 0 to $x_i$ (which corresponds to choose the  particular solution vanishing at the origin):
$$ u_h(x_i)\simeq\int_0^{x_i} (u_h)_x\ dx = \sum_{k=1}^{i} \int_{x_{k-1}}^{x_{k}} z_{k}= h \sum_{k=1}^{i} z_{k}. $$
In higher dimension a different strategy can be used:  plugging $\nabla\psi_h=v_h\nabla u_h$  into (\ref{discr_neum}), $u_h$ becomes the solution of the discrete variational Neumann problem (with $v_h$ given)
\begin{equation}
u_h\in V_h, \quad \int_{\Omega_h} v_h \nabla u_h \cdot \nabla \phi \ dx = \int_{\Omega_h} g_h \phi \ dx\ , \quad \forall \phi\in V_h .
\end{equation}

\mascotsection{Computation of growing profiles}

In  this section we want to show that the previously characterized similarity solutions  asymptotically arise as profiles of the growing heaps in the dynamic process of filling the silo. In order to do that we implemented a numerical scheme for the complete system (\ref{hk}), adapting to this case the finite difference scheme  used in \cite{FFV06} for the growing sandpiles. 

In 1D, if $\Omega=(0,L)$ and $h=L/(N-1)$ denotes the space discretization step, a uniform mesh is described by the nodes $x_i=(i-1)h$, for $i=1,..,N$. If $\Delta t$ is the time step, our explicit scheme reads for the internal nodes $i=2,..,N-1$ as
\begin{equation}
\left\{\begin{array}{ll}
v^{n+1}_i&=v^n_i+\Delta t(G_i-(1-|Du_i^n|)v^n_i+f_i)\\ \\
u^{n+1}_i&=u^n_i+\Delta t(1-|Du_i^n|)v^n_i
 \end{array}\right.
\end{equation}
where  $u^n_i, v^n_i,G^n_i$  denote  respectively the approximate values at time $n\Delta t$ in $x_i$ of the solutions $u,v$ and of the upwind flux derivative $(vu_x)_x$  in the direction determined by the sign of $Du_i$, that is 
$$G^n_i=\left\{\begin{array}{ll}
(v^{n}_{i+1}Du^n_{i+1}-v^n_iDu^n_i)/h  \quad \hbox{ if } Du^n_i>0\\ \\
(v^{n}_{i}Du^n_{i}-v^n_{i-1}Du^n_{i-1})/h  \quad \hbox{ if } Du^n_i<0\ 
 \end{array}\right. $$      
(in each node the spatial derivative $Du_i$ is defined as the term of maximal absolute value between the backward and the forward first differences). To complete the scheme we added initial conditions ($u^0_i=v^0_i=0 \ \forall i$)
and  boundary terms induced by the Neumann condition on $u$.

The extension of this approach to the 2D case is straightforward if we restrict ourselves to square or rectangular cross-sections for the silo. It is enough to decompose the flux term as
$\nabla\cdot(v\nabla u)= (vu_x)_x + (vu_y)_y $,
and to repeat the 1D approach in each direction. 

in order to study the asymptotic behavior of the growing heaps, in the numerical tests the scheme was stopped when the relative growth per iteration of the standing layer resulted approximately the same at each node, revealing the emergence of a similarity profile. Such profile was then translated to the base of the silo and compared with the computed similarity solution.

\begin{figure}[h] 
\center{\fbox{\psfig{file=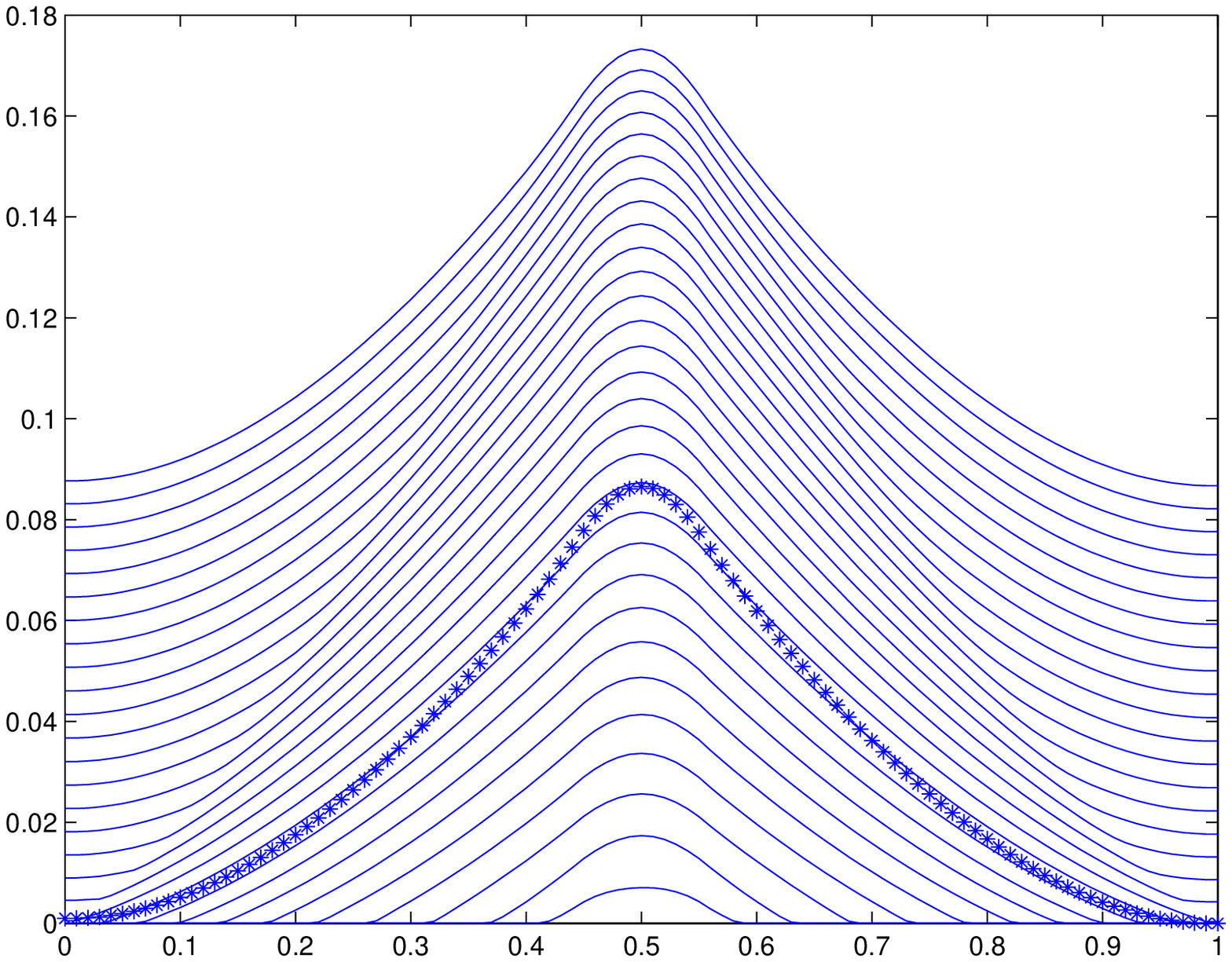,height=1.1in}}\hspace{.005in}
\fbox{\psfig{file=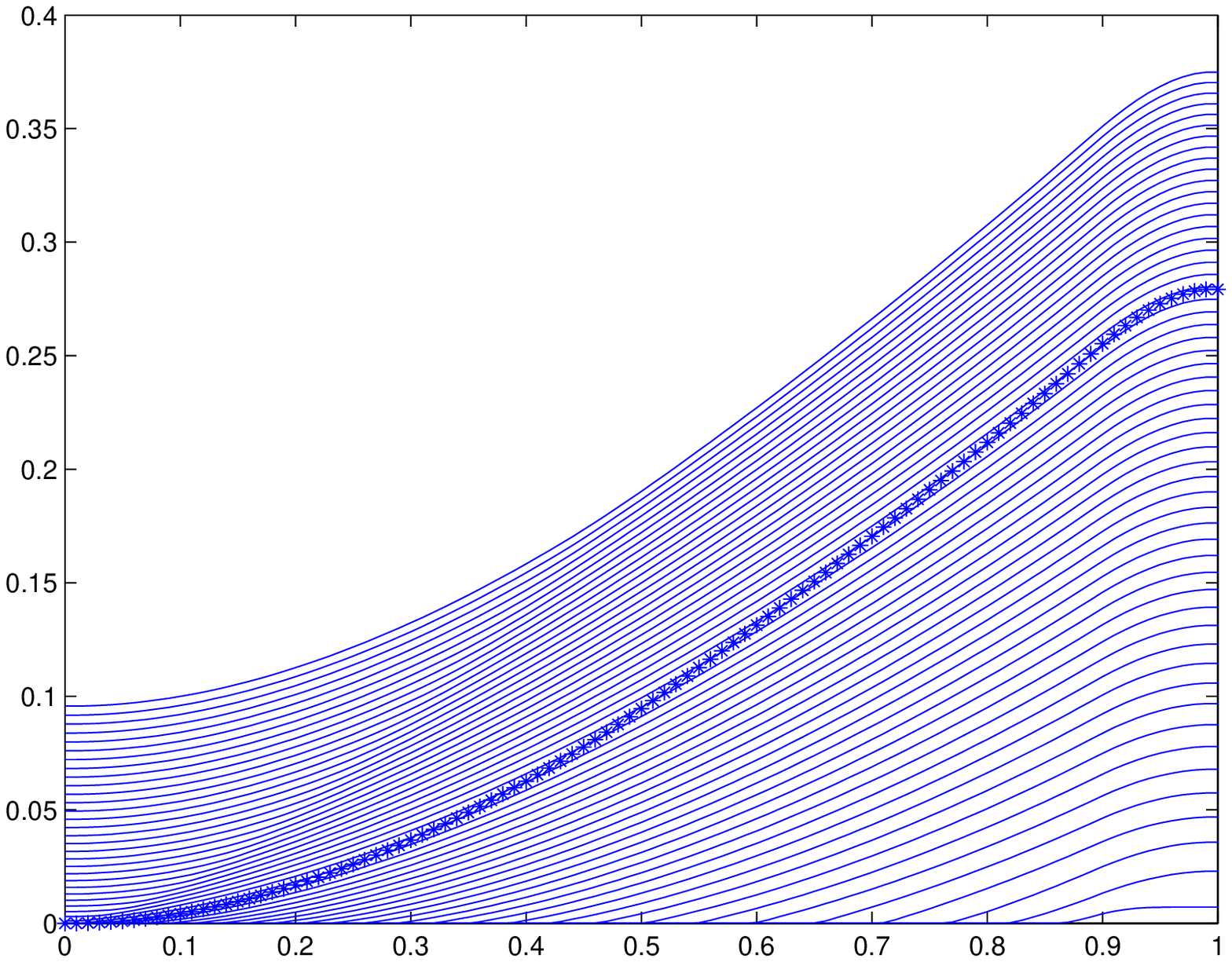,height=1.1in}}\hspace{.005in}
\fbox{\psfig{file=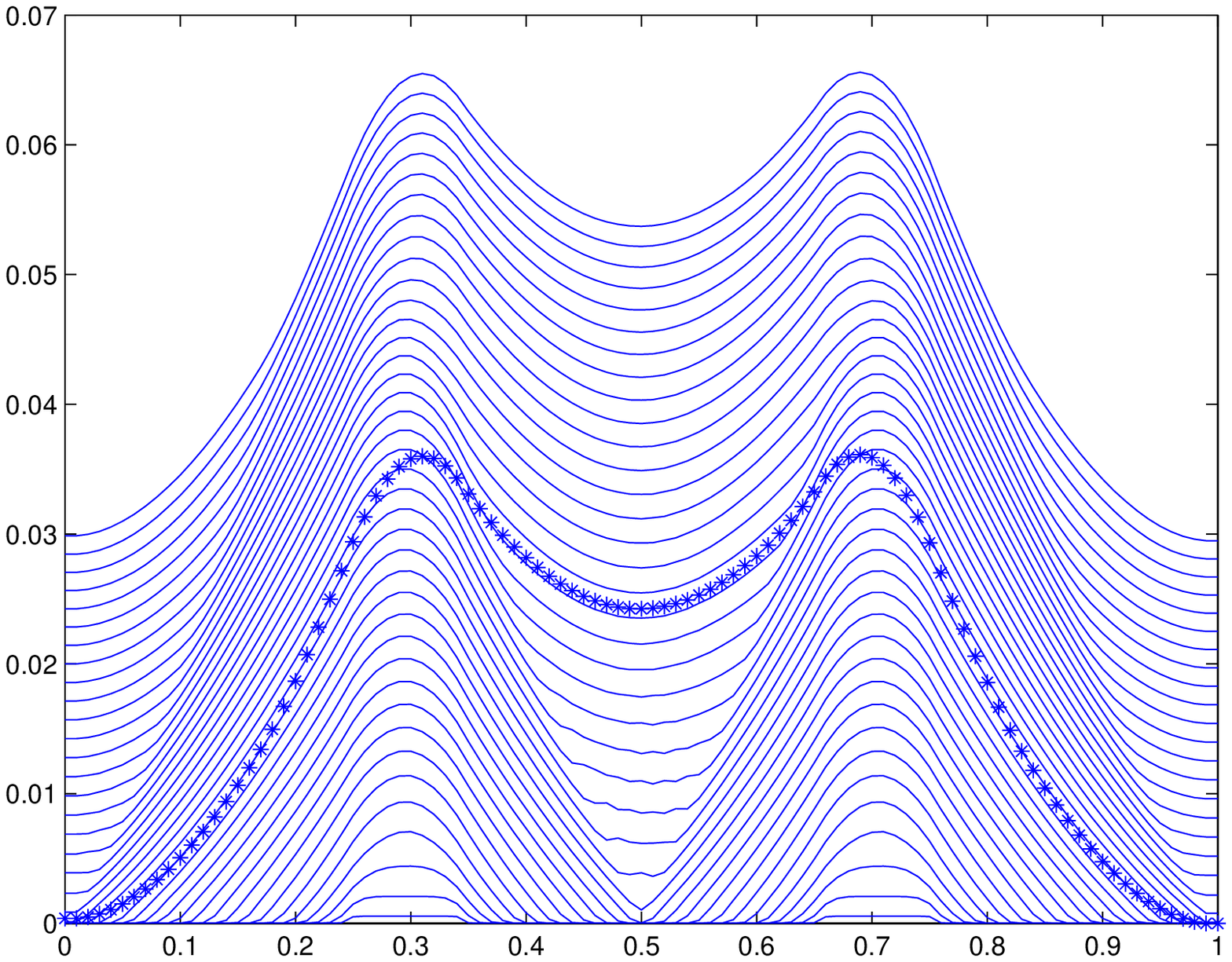,height=1.1in}}}
\caption{Growing heaps and similarity profiles in 1D silos with different source supports $S_f\subseteq [0,1]$: a)[0.45,0.55]; b)[0.9,1]; c)[0.25,0.35]$\cup$[0.65,0.75].}
\label{grow1d}
\end{figure}

\mascotsection{Numerical experiments}
For the 1D cross-section case we assumed $\Omega=(0,1)$ and  tested different choices of the source term $f$. In each example, for a given uniform partition (of step $h$) of $\Omega$, we compared the exact similarity solutions given by formulas (\ref{sim1}) (the couple $(U, V)$, with  $\min_\Omega U=0$), the discrete similarity solutions computed by the FE approach of Section 3 ($(u^e,v^e)$, with $\min_i u^e_i=0$), and  the stabilized profiles determined by the FD scheme for the evolutive problem described in Section 4 ($(u^d,v^d)$, with the first one shifted towards zero, that is $u^d=u^n-\min_i (u_i^n)$, where $u^n$ is the iterate selected by the stopping criterion). 
Figure \ref{grow1d}  shows three examples of growing heaps according to different source supports $S_f$ (centered symmetric, close to the boundary, disconnected). In all the cases it can be seen the formation of a similarity solution (in dot lines). 

\begin{figure} 
\center{\fbox{\psfig{file=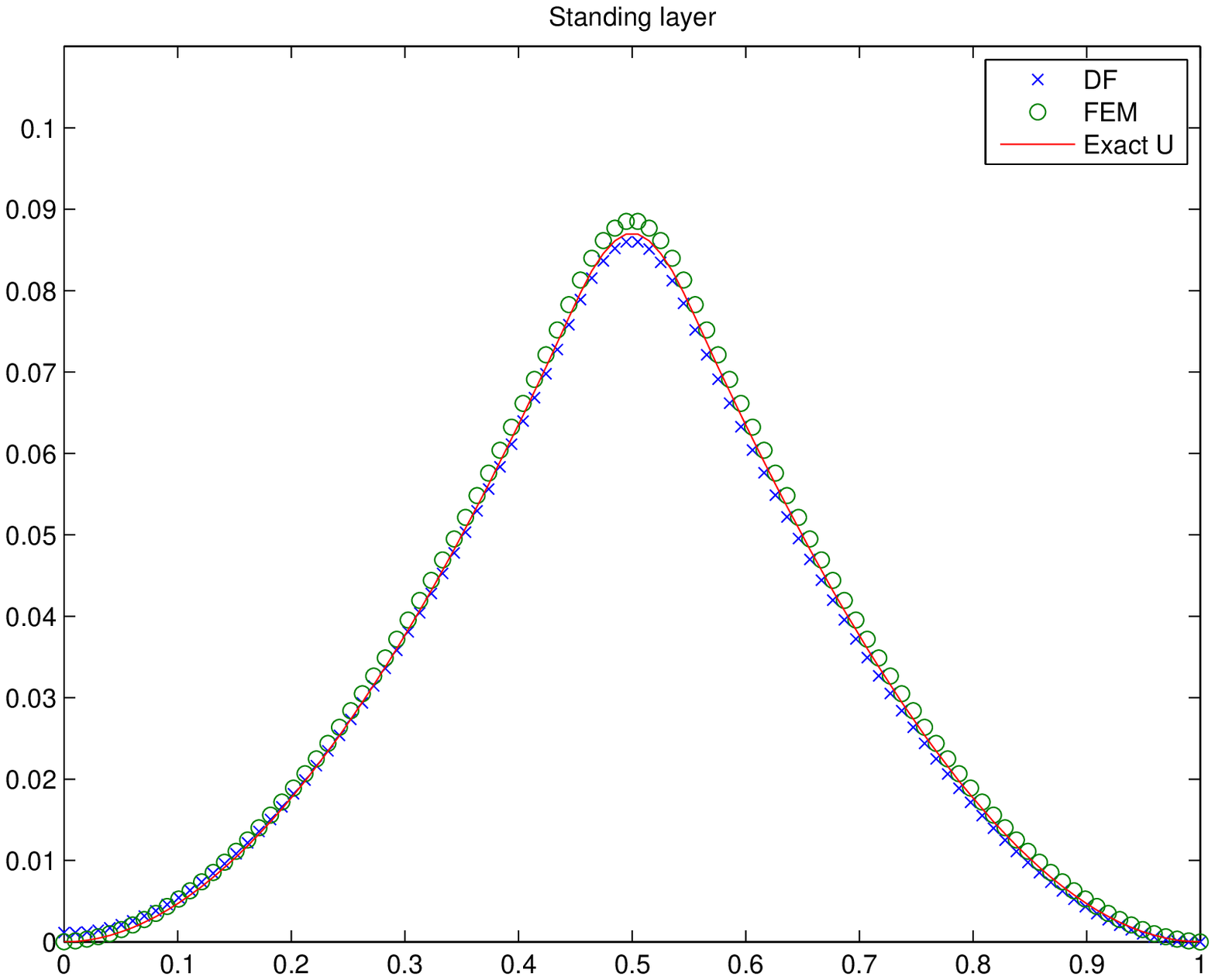,height=1.in, width=1.5in}}\hspace{.1in}
\fbox{\psfig{file=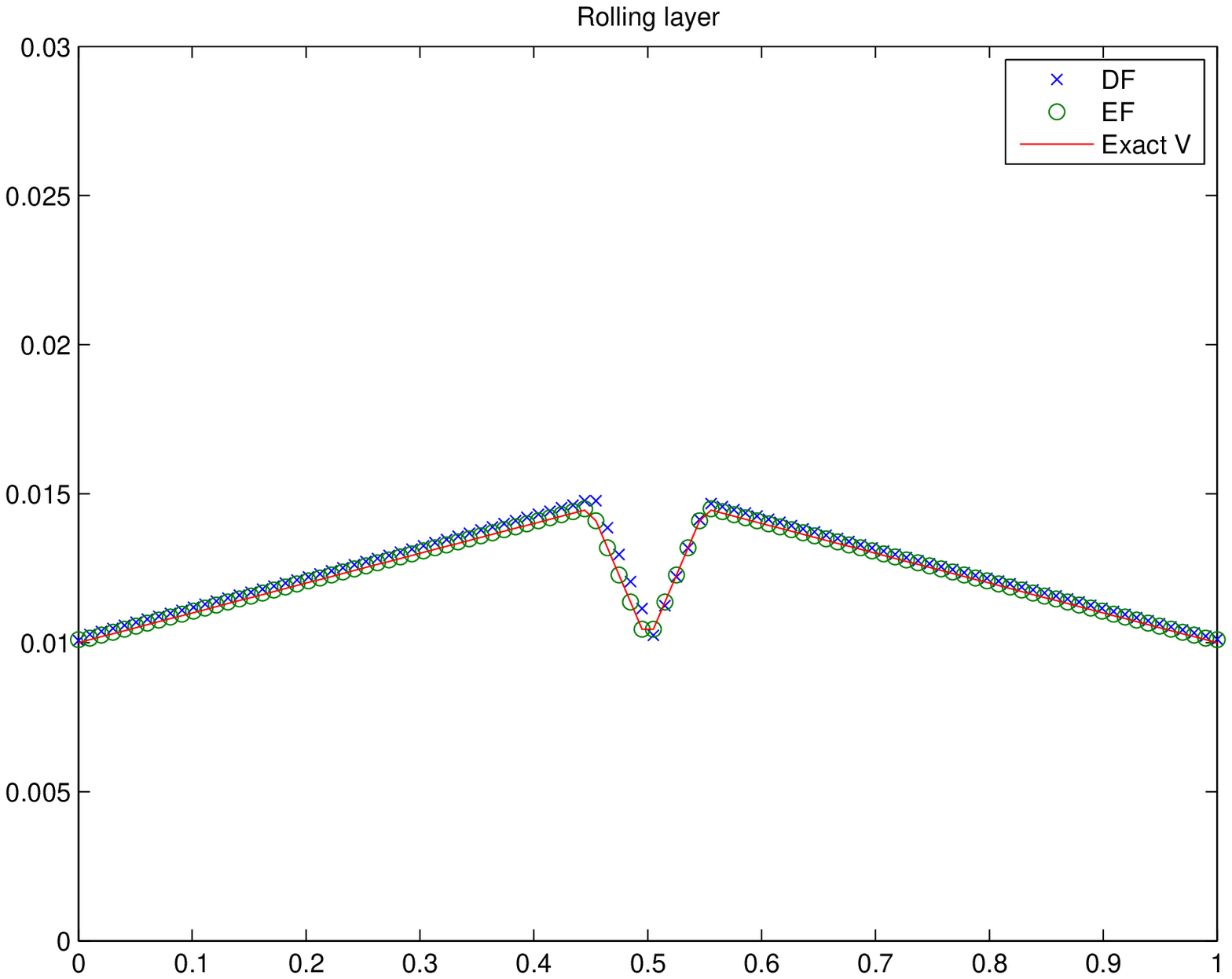,height=1.in, width=1.5in}}}
\caption{Standing and rolling layers for 1D silo with small central  source.}
\label{sol1dcentral}
\end{figure}

We found first order convergence of the approximate similarity solutions to the exact ones for the FE method in uniform norm, and approximately the same order for the asymptotic convergence of the growing profiles computed by the FD scheme to the quasi-stationary ones. In Table 1 we report the values found for the symmetric centered support case of Figure \ref{sol1dcentral}. Other tests gave similar results. Note that  the stabilized rolling layer is everywhere positive, with a small depression in the central region corresponding to the source support. When its length tends to zero one recovers the situation of the point source of Example~1, that is $u$ assumes the known logarithmic profile and the depression region of $v$  disappears.
\begin{table}
\centering
\begin{tabular}{|c|c|c|c|c|}
\hline
h      & \multicolumn{2}{c|}{Standing layer} & \multicolumn{2}{c|}{Rolling layer} \\ \hline
       & $\|U-u^e\|_\infty$             & $\|U-u^d\|_\infty$              & $\|V-v^e\|_\infty$            & $\|V-v^d\|_\infty$            \\ \hline
0.01   &   $3.07\times 10^{-3}$             &    $ 4.97 \times 10^{-3 } $          &    $4.6 \times 10^{-4 }  $          &   $1.62 \times 10^{-3 }  $          \\ \hline
0.005  &   $1.55\times 10^{-3 } $         &      $ 2.66\times 10^{-3 }   $       &      $ 2.3\times 10^{-4 } $         &    $  8.4\times 10^{-4 }    $      \\ \hline
0.0025 &   $7.8\times 10^{-4 } $             &    $ 1.44  \times 10^{-3 } $         &    $  1.2\times 10^{-4 } $          &  $  5.3 \times 10^{-4 }    $       \\ \hline
0.001  &   $ 3.1\times 10^{-4 } $            &    $ 4.7 \times 10^{-4 }  $         &    $  5\times 10^{-5 }  $         & $   2\times 10^{-4 }  $          \\ \hline
\end{tabular}
\caption{$L^\infty$ errors for EF and DF schemes in the test case of Figure \ref{sol1dcentral}.}
\end{table} 

In the more realistic case of a spatial silo, that is when the cross-section $\Omega$ is a 2D domain,  the similarity solutions can only be approximated, so that we just estimated the quantities $ \| u^e-u^d\|_\infty$ and $\| v^e-v^d\|_\infty$, where as before $(u^e,v^e)$ denote the FE solutions of the stationary model and $(u^d,v^d)$ the FD solutions of the evolutive model after the stopping criterion applies, computed over the same mesh introduced in $\Omega$. We restricted our tests to the case of a rectangular domain $\Omega$ decomposed through a uniform mesh, in order to use the same set of nodes for the two schemes. If for example $\Omega=(0,1)\times(0,1)$, a mesh with $N\times N$ equispaced nodes (with $h=\Delta x=\Delta y=1/(N-1))$ for the FD scheme can be used as well as a base for a uniform Courant FE triangulation over $\Omega$. 

\begin{figure}[hp] 
\center{\fbox{\psfig{file=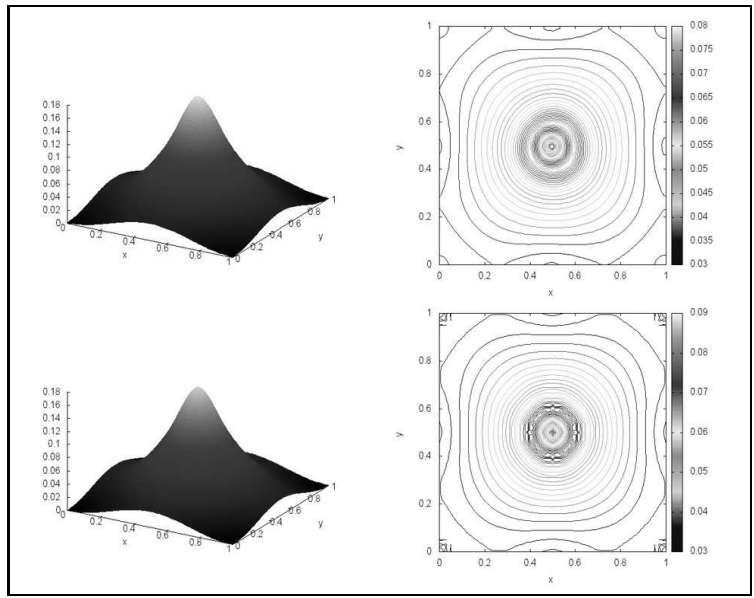,height=3.in, width=3.8in}}}
\caption{Similarity solutions (above) and asymptotic profiles (below) in a square silo: source supported in a central small ball.}
\label{cen2d}
\end{figure}

\begin{figure} 
\center{\fbox{\psfig{file=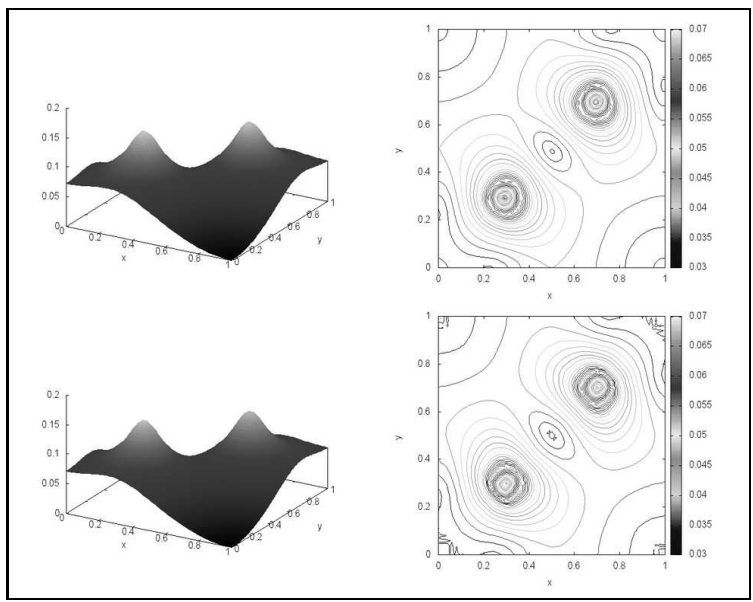,height=3.in, width=3.8in}}}
\caption{Similarity solutions (above) and asymptotic profiles (below) in a square silo: source supported in the union of two disconnected balls.}
\label{scon2d}
\end{figure}

\begin{figure} 
\center{\fbox{\psfig{file=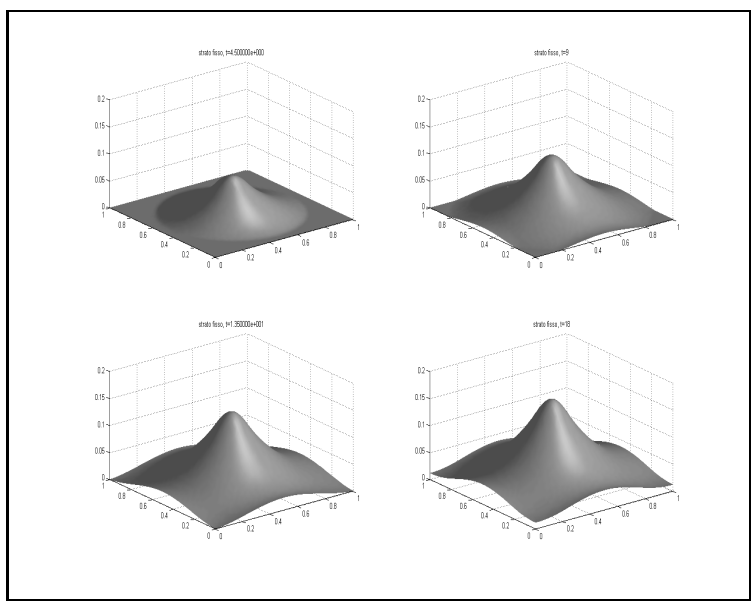,height=3.in, width=3.8in}}}
\caption{Growing heap in a square silo: source supported in a central ball.}
\label{riemp}
\end{figure}

The experiments gave results similar to those of the 1D case. Figures \ref{cen2d} and \ref{scon2d} illustrate the results corresponding to a source supported in a small ball in the center of the silo or in the union of two disconnected balls, showing the profiles of the standing layers and the level lines of the rolling layers.  The correspondence of the similarity solutions (above) to the evolving profiles (below) appears evident. Figure \ref{riemp} shows the growing heap in the square silo for the first example at four successive time steps.

\end{document}